\documentclass[12pt]{article}
\usepackage{amssymb, amsmath, amsthm}
\def\Rep{{\rm Rep\,}}
\def\Hom{{\rm Hom}}
\def\Ob{{\rm Ob\,}}
\def\Im{{\rm Im\,}}
\def\rank {{\rm rank\,}}
\def\ind{{\rm ind\,}}
\def\mod{{\rm mod\,}}
\def\aff{{\rm aff\,}}
\def\Ker{{\rm Ker\,}}
\def\El{{\rm El\,}}
\def\Mor{{\rm Mor\,}}
\def\Sets{{\mathit{Sets}}}

\renewenvironment{proof}%
    {\par\noindent\textit{Proof.}\ }{\par\medskip}

\theoremstyle{plain}
\newtheorem{theorem}{Theorem}
\newtheorem{lemma}{Lemma}
\newtheorem{corollary}{Corollary}
\newtheorem{proposition}{Proposition}

\theoremstyle{definition}

\newtheorem{remark}{Remark}
\newtheorem{example}{Example}

\newcounter{pgraphno}
\newenvironment{pgraph}%
   {\begin{list}{\bf{\arabic{pgraphno}}. }{\usecounter{pgraphno}%
    \setlength{\labelsep}{0pt}\setlength{\leftmargin}{0pt}%
    \setlength{\labelwidth}{0pt}%
    \setlength{\parsep}{0pt}%
    \setlength{\listparindent}{\parindent}}}%
   {\end{list}}

\newcommand{\citem}{%
    \setcounter{lemma}{\thepgraphno}%
    \setcounter{theorem}{\thepgraphno}%
    \setcounter{proposition}{\thepgraphno}%
    \setcounter{corollary}{\thepgraphno}%
    \setcounter{remark}{\thepgraphno}%
\item}

\title{Representations of marked quivers}
\author{A.V. Roiter\thanks{e-mail: roiter@imath.kiev.ua} \\
\small Institute of Mathematics of National Academy of Sciences of
Ukraine,\\ \small 3 Tereshchenkivska Street,  01601 Kyiv-4,
Ukraine}
\date{}

\begin{document}

\maketitle

We introduce a generalization of representations
of quivers~[1] that contains also
representations of posets~[2], vectorspace
problems and other matrix problems. Many
examples, some of which are given below, show
that the language of marked quivers is rather
convenient and that the notion of their
representations can probably be considered as
basic for a general theory of matrix problems.
Some historical comments from my personal point of view
are added in the Appendix.
I am very grateful to the Referee for the numerous remarks and advices.

We compose maps on the right. Thus  we denote the composition of
$\alpha: A \longrightarrow B$ and $\beta: B \longrightarrow C$ by
$\alpha \beta$.

\begin{pgraph}
\citem Let $Q$ be a quiver, $Q_{v}$ the set of vertices (points)
of $Q$, $Q_a$ the set of arrows. We recall that a representation
$U$ of $Q$ in a category $K$ attaches to any $a \in Q_v$ the
object $U(a) \in K$ and to any $x \stackrel {\alpha}
{\longrightarrow} y$ the morphism $U(\alpha)\in K(U(x), U(y))$.
Representations of $Q$ in $K$ form the category $\Rep Q = \Rep
(Q,K)$, and a morphism $f: U\longrightarrow W$ between the
representations $U$ and $W$ consists of morphisms $f(x): U(x)
\longrightarrow W (x)$, $x \in Q_v$ such that $f(x) W(\alpha) =
U(\alpha) f(y)$ for each arrow $x
\stackrel{\alpha}{\longrightarrow} y$. Usually $K$ is the category
of modules over a field (or a commutative ring) $k$ [1, 3].

We will denote by $\Delta $ or $\Delta_{\alpha}$ the quiver $x
\stackrel{\alpha}{\longrightarrow} y$, ($x \neq y$), $\Delta_v
=\left\{x,y\right\}$, $\Delta_{a} = \left\{\alpha \right\}$. The
category $\Rep (\Delta ,K)$ is called in~[4] the
{\it category of morphisms} of a category $K$.

We will say that a quiver $Q$ is {\it marked} if to each
$x \in Q_{v}$ is attached a category $K_{x}$, and to each arrow $x
\stackrel{\alpha}{\longrightarrow} y$
 is attached a functor  $\Phi_{\alpha }: K^{\circ}_{x} \times K_{y}
 \longrightarrow  \Sets$, (i.e., a bifunctor contravariant in $K_{x}$
 and covariant in $K_{y}$ and taking values in the category of sets).
Then $M = \left\{K_{x},\Phi_{\alpha}\right\}$ is a
{\it marking} of the quiver $Q$.

A representation $U$ of a marked quiver $Q_{M}$ attaches to each $a
\in Q_{v}$ the object $U(a) \in K_{a}$ and to each $x
\stackrel{\alpha}{\longrightarrow} y$ the element $U(\alpha ) \in
\Phi_{\alpha}(U(x),U(y))$. A morphism $f:U \longrightarrow W$
consists of morphisms $f(x) \in K_{x}\left(U(x), W(x)\right)$, $x\in
Q_{v}$ such that $W(\alpha )\Phi_{\alpha }\left(f(x), 1_{y}\right)
= U(\alpha ) \Phi \left(1_{x},f(y)\right)$ for each arrow $x
\stackrel{\alpha}{\longrightarrow} y$.

Multiplication of morphisms (in both $\Rep Q$ and $\Rep Q_M$) is
defined in the natural way.
\begin{remark} For any functor $\Phi : K^{\circ} \times L \longrightarrow
{\Sets}$, $t \in \Phi(A,B)$, $\alpha \in K(A',A)$ and $\beta\in
L(B,B')$, where $A,A' \in \Ob K$ and $B,B' \in \Ob L$, we can
agree  to define the compositions $\alpha t \in \Phi(A',B)$ and
$t\beta \in \Phi(A,B')$ by putting $\alpha t=t\Phi(\alpha, 1_B)$
and $t \beta= t\Phi(1_A,\beta)$. In fact by such agreement we
attach (in the additive case) to the functor $\Phi$
 a {\it bimodule} ${_K}{\Phi}_L$
over the categories $K$ and~$L$ (compare [3], 2.2). With this
agreement the formulas in the definitions of categories of
representations of unmarked and marked quivers coincide, and the
second of these definitions repeats the first one word-for-word.
\end{remark}

In order to obtain representations of  the ``usual'' (unmarked)
quiver $Q$ in a category $K$ we should of course put $K_{x} = K$
for all $x \in Q_{v}$, and $\Phi_{\alpha}  = \Hom_{K}$ for all
$\alpha \in Q_a$.

We will say that a marking $M$ is a {\it point-marking} (in
a category $K$) if there exist  a category $K$ and functors
$\Phi_{x}: K_x\longrightarrow K (x\in Q_v)$ such that
$\Phi_{\alpha}(A,B)=K(\Phi_{x}(A),\Phi_{y}(B))$ and $t
\Phi_{\alpha}(p,q)=\Phi_{x}(p)\> t\> \Phi_{y}(q) \in
\Phi_\alpha(A',B')$ for any $x \stackrel{\alpha}{\longrightarrow}
y$ (may be $x=y$) and $t\in \Phi_\alpha(A,B)$. Here $A,A' \in \Ob
K_x$, $B,B' \in \Ob K_y$, $p\in K_x(A,A')$ and $q\in K_y(B,B')$.

\citem Let $Q_M$ be a marked quiver and $Q^i$ a subquiver of $Q$.
The marked quiver $Q^i_M$ is defined naturally by restriction of
the marking $M$ to $Q^i$. For $x\in Q_v$, we define a functor
$Y_x : \Rep Q_M \to K_x$, $U\mapsto U(x)$, and
for $x\in Q_v^i$ a functor  $Y^i_x : \Rep Q^i_M \to K_x$.

Let  $\omega$ be a set, $\left\{Q^{i} | i\in \omega \right\}$ be a
family of subquivers of $Q$ such that $Q_{v} = \bigsqcup_{i \in
\omega } Q^{i}_{v}, \quad  Q_{a}^{i} \bigcap Q^{j}_{a} ={\varnothing}$, if
$ \left\{i, j \right\} \subset \omega$, $i\ne j $. We introduce a
quiver $\Omega$ with vertex set $\Omega_v = \omega$ and arrow set
$\Omega_a$ in one-one correspondence with $Q_a\setminus
\bigcup_{i\in \omega} Q^i_a$ and defined as follows. Corresponding
to each arrow $x\stackrel{\alpha}{\to} y \; \in \; Q_a\setminus
\bigcup Q^i_v$ with $x\in Q^i_v$ and $y\in Q^j_v$ there is an
arrow $i\stackrel{\overline{\alpha}}{\to} j\;\in\;\Omega_a$; if
$i$ and $j$ coincide, then $\overline{\alpha}$ is a loop.

Let $M$ be a marking of $Q$. We define the marking $\overline{M}$
of $\Omega$,  by putting $\overline{K_{i}} = \Rep Q^{i}_{M}$,
$\Phi_{\overline{\alpha}} = \left(Y^{i}_{x} \times Y^{j}_{y}
\right) \Phi_\alpha $, where $x \stackrel{\alpha}{\longrightarrow}
y$, $x \in Q^{i}_v$, $y \in Q^{j}_v$. From these definitions
directly follows
\begin{proposition}
$\Rep Q_{M} \simeq \Rep \Omega_{\overline{M}}$.
\end{proposition}
For $z\in Q_v$, let $\delta(z)$ denote the number of arrows
incident on it.

Now consider the following example. Suppose that there is a vertex
$z\in Q_v$ with $\delta(z) = 1$ and that $\beta$ is the unique
arrow incident on it. Let $w\neq z$ be the second vertex of
$\beta$ (either $z \stackrel{\beta}{\to} w$ or $w
\stackrel{\beta}{\to} z$). We denote by $Q'$ the subquiver of $Q$
such that $Q'_v=Q_v\setminus\{z\}$ and $Q'_a =
Q'_a\setminus\{\beta\}$.

For any $x_i \in Q_v$ we define the subquiver $Q^i \subset Q$ as
follows. If $x_i \neq w$, then $Q^{i}_{v}=\left\{x_i\right\}$,
$Q^i_{a}={\varnothing};$ and if $x_i = w$, then $Q^i=\Delta_\beta$. We
identify the quiver $\Omega$ with $Q'$. Let $\overline{M}$ be the
marking induced on $Q'$, so that $\overline{K_x} = K_x$ if $x \neq
w$, $\overline{K_w} = \Rep(\Delta_\beta)_M$ and the functors are
defined naturally.
\begin{corollary}
$\Rep {Q_M} \simeq \Rep Q'_{\overline M}$.
\end{corollary}
\begin{remark}
If $M$ is a point-marking then $\overline M$ in the corollary (but
not in the proposition!) is also a point-marking.
\end{remark}
\citem Let $k$ be a fixed algebraically closed field. By $\mod k$
is denoted the category of finite dimensional vector spaces over
$k$. Unless otherwise stated all categories will be assumed to be
{\it $k$-categories} (i.e. the sets $\Hom (A,B)$ are endowed
with a $k$-module structure such that the composition of maps is
$k$-bilinear), all subcategories to be
{\it $k$-subcategories}, and all functors to be
{\it $k$-functors}~[3]. In particular, the functors $\Phi
_\alpha$ of the definition of a marked quiver are $k$-functors
$K^{\circ}_x \times K_y \longrightarrow \mod k$ (which, of course,
may be also viewed as functors to {\it Sets}).

An additive $k$-category in which every idempotent has a kernel is
an {\it aggregate} [3]. A category $S$ is a
{\it spectroid} if its objects are indecomposable (i.e., for
any $A \in \Ob S$, $S(A,A)$ contains exactly two idempotents $1_A
\neq O_A$) and pairwise non-isomorphic~[3]. A subspectroid of
$\mod k$ (i.e. a $k$-subcategory of $\mod k$, which is a
spectroid) is a {\it vectroid}~[5]. Spectroids are not
additive, but any spectroid $S$ generates in a natural way its
``additive closure'', the aggregate $\oplus S$~[3] (whose objects
are the sequences $(X_1,\dots,X_m)$, $X_i \in \Ob S$). On the other
hand, any aggregate $\mathfrak{A}$ is defined by its spectroid
$S(\mathfrak{A})$, (such that $\mathfrak{A} \simeq \oplus
S(\mathfrak{A})$).

If $V$ is a vectroid, the aggregate $\oplus V$ is naturally
embedded in $\mod k$.

For a functor $\Phi: \mathfrak{A} \longrightarrow \mathfrak{B},
\Ker \Phi $ is the set (an ideal of $\mathfrak{A}$) of morphisms
$\gamma$ of $\mathfrak{A}$ such that $\Phi(\gamma)= 0;$ if $A \in
\Ob \mathfrak{A}$ then $A \in \Ker\Phi$,  means that $1_A \in
\Ker\Phi$. The functor $\Phi$ is {\it faithful} if $\Ker
\Phi = \{0\}$, then $\mathfrak{A} \simeq \Im \Phi$.

We will say that $Q$ is {\it $k$-marked} if it is point-marked in
$\mod k$,  all $K_x$ are aggregates,  and the functors $\Phi _x$
$(x\in Q_v)$ are faithful. So in this case we may consider $K_x$
$(x \in Q_v)$ as subaggregates of $\mod k$.

Choosing a representative in each indecomposable isoclass of $K_x$
$(x \in Q_v)$, we obtain a vectroid $V_x$. So a $k$-marked quiver
is defined by the collection of vectroids $V_x$, $x \in Q_v$ $(K_x =
\oplus V_x)$.

Unless otherwise stated we will consider only {\it finite}
$k$-markings, i.e., we will assume, for a $k$-marked quiver $Q_M$,
that $\left|Q_v\right|< \infty$, $0< \left|Q_a\right|<\infty$,
$\left|\Ob V_x \right|< \infty$, $(x \in Q_v)$. We will also consider
only connected quivers.

For a category $C$, we denote by $\ind C$ the set of isoclasses of
indecomposables in $\Ob C$; we abbreviate $\ind\Rep Q_M$ to $\ind
Q_M$. We say that $Q_M$ has {\it finite type} if
$\left|\ind Q_M \right| < \infty$.

For a fixed $x \in Q_v$, $T \in {\ind K_x}$, we construct $ U_T \in
\ind Q_M : U_T (x) = T$, $U_T (y) = 0$ for $x\neq y$,
$U_{T}(\alpha) = 0$ for any $\alpha \in Q_a$. So we may consider
$\bigsqcup_{x \in Q_v} \ind K_x \subset \ind Q_M$.

If $M$ is a point-marking of $Q$ in $\mod k$ such that the $K_x$
are aggregates but the functors $\Phi_x$ are not faithful,  we may
define a $k$-marking $\tilde{M}$ (of $Q$) by putting $\tilde{K_x}
\simeq K_x \diagup \Ker \Phi_x$ and taking $\tilde{\Phi}_x$ to be
the embedding of $\tilde{K}_x = \Im \Phi_x$ into $\mod K$. The
next statement is obvious
\begin{lemma} There is a natural bijection between $\ind Q_M$ and
$\ind \tilde{Q}_{\tilde M} \sqcup \bigsqcup_{x \in Q_v} (\ind K_x
\cap \Ker \Phi_x)$.
\end{lemma}
\citem
If $\mathfrak{A}$
and $\mathfrak{B}$ are two categories (in particular, vectroids)
then $\mathfrak{A} \sqcup \mathfrak{B}$ is the category such that
$\Ob (\mathfrak{A} \sqcup \mathfrak{B} ) = \Ob \mathfrak{A} \sqcup
\Ob \mathfrak{B}$, $(\mathfrak{A} \sqcup \mathfrak{B} )(A_1,A_2) =
\mathfrak{A} (A_1,A_2)$, $(\mathfrak{A} \sqcup
\mathfrak{B})(B_1,B_2) = \mathfrak{B} (B_1,B_2)$,
$(\mathfrak{A} \sqcup \mathfrak{B}) (A_1,B_1) = \left\{0\right\}$,
$(\mathfrak{A} \sqcup \mathfrak{B}) (B_1,A_1) = \left\{0\right\}$
for $A_1,A_2 \in \Ob \mathfrak{A}$, $B_1, B_2 \in \mathfrak{B}$.
\begin{example}  Let $P$ be a poset considered as a
category ($\Ob P = P$,  $\left|\Hom (u,w)\right| = 1$, if $u \leq w$,
$\Hom (u, w) = {\varnothing}$ otherwise). Let $ kP $ be its {\it $k$-
linearization}~[3], which may be viewed as a vectroid (to each $p
\in P$ is associated a one-dimensional vector space $kp \in \mod
k$, $\dim \Hom_{kP}(a,b) = 1$ if $a \leq b$, $\Hom_{kP}(a,b) =
\left\{0\right\}$ otherwise).

We define two vectroids, $k_n$ and $k^n$, as follows. The vectroid
$k_n$, which we call {\it linear}, is the linearization of a
linearly ordered set $P$ with $|P|=n$. On the other hand, the
vectroid $k^n$ has one object $A$ and morphism space $\Hom(A,A) =
k\langle{r}\rangle$ with $r^n=0$. We take $A$ to be an
$n$-dimensional $k$-space with basis $\{a_1,\ldots,a_n\}$ and
define an action of $r$ by $a_ir=a_{i+1}$, $a_nr=0$. Observe that
$k_1=k^1=k$. We say that the vertex $x$ is an {\it unmarked
point} if $V_x = k$ (so that $K_x \simeq \mod k$).

If $\Delta = x \stackrel{\alpha}\longrightarrow y$,  $V_x=k$,  $V_y =
kP$, then the category $\Rep \Delta_M$ coincides with the category
of representations of the poset $P$ in the sense of~[2]
($\bigoplus V_y \subset \mod K$, $\Phi_\alpha = \Hom_k$, elements of
$\Hom (X,Y)$, where $X\in \mod k$, $Y\in \bigoplus V_y$, may be
considered as matrices divided into vertical strips indexed by
elements of $P$).
\end{example}
\begin{example}
If  $Q=\Delta$, $V_x = k$ and $V_y$ is an arbitrary vectroid $V$,
then $\Rep \Delta_{M_V}$ coincides with the category $\Rep V$ for
the vectroid $V$ in the sense of~[5] (= categorical matrix
problem~[6]= vectorspace problem~[7]), and we identify these two
categories. So $\ind\, V$~[5] is $\ind \, \Delta_{M_V}$, $V$ is
{\it tame} if $\Delta_{M_V}$ is (see Section~5) and so on. If $M$ is
an arbitrary $k$-marking of $\Delta$ then $\Delta _M$ coincides
with the category of representations of the pair $(V_x,V_y)$~[8].
If $V_y=k$ and $V_x=V$, then $\Rep \Delta _M\simeq \Rep
V^{\circ}$, where $V^{\circ}$ is the vectroid {\it opposite} to
$V$ (i.e. there are bijections $*$ between $\Ob V$ and $\Ob V^*$,
between $ A$ and $A^*$ for $A\in V$, between $ \Hom_V (A,B)$ and
$\Hom_{V^*}(B^*,A^*)$ such that if $a\varphi = b$, then
$b^*{\varphi ^*}=a^*$ for $a\in A$, $b\in B$, $\varphi \in
\Hom_V(A,B)$).
\end{example}
\begin{example}
We recall [3] that a right (finite dimensional) module
$\mathfrak{M} $ over a $k$-category $\mathfrak{A} $ consists of a
finite dimensional vector space $\mathfrak{M} (X)$ (over $k$) for
each object $X\in \mathfrak{A}, $ and of maps $\mathfrak{M}
(X)\times \mathfrak{A} (X,Y)\longrightarrow \mathfrak{M}
(Y)$, $(m,f)\longrightarrow mf$ that satisfy the usual axioms
($m1_x=m$, $m(fg)=(mf)g (m_1+m_2)f=m_{1}f+m_{2}f$ and
$m(f_{1}+f_{2})=mf_{1}+mf_{2}$).

  The notions of {\it submodule}, {\it factormodule}, and
  so on, are defined in the natural way. A functor $\Phi _\mathfrak{M} :
   \mathfrak{A} \longrightarrow \mod k$ arises from the module
$\mathfrak{M} $. Then the category $\mathfrak{M}^k$ defined in [3]
coincides with $\Rep\Delta _M$ where $M$ is a point-marking in
$\mod k$ (not a $k$-marking) of $\Delta :x$ is unmarked,
$K_y=\mathfrak{A}$, $\Phi _y=\Phi _\mathfrak{M} $. Lemma~3 holds.
\end{example}
\begin{example}
We call a family $\mathfrak{M} _1,\dots,\mathfrak{M} _p$ of
modules over an aggregate $\mathfrak{A} $ {\it a bunch of
modules}. The category of representations of a bunch of modules is
$\Rep Q_M$, where $Q_v = \{x, y_1,\dots,y_p\}$, $Q_a
=\{\alpha_1,\dots,\alpha_p\}$, $x \stackrel {\alpha _i}
{\longrightarrow} y_i$, $K_x=\oplus k$, $K_{y_i}=\mathfrak{A}$
$(1\leq i\leq p)$, $\Phi _{\alpha_i}:(\oplus k)^{\circ}\times
K_{y_i}\longrightarrow \mod k$ are induced by $\Phi
_{\mathfrak{M}_i}:\mathfrak{A} \longrightarrow \mod k$ $(i =
1,\dots, p)$~[8]. In this case $M$ is not a point-marking.
\end{example}
\begin{example}
Let $Q$ be $c\longleftarrow d\longrightarrow b$. A representation
of the (non-marked) $Q$ may be viewed as a matrix over $k$ divided
into two vertical strips.
\begin{center}
\begin{tabular}{c|c|c|}
\multicolumn{3}{c}{\ \ \ \ \ $c$ \ \ $b$} \\
\cline{2-3}
$d$ & \ \ \ & \ \ \ \\
  \cline{2-3}
\end{tabular}
\end{center}
Two representations  are equivalent if one can
be obtained from the other by means of any row
transformation and any column transformations
inside of each strip.

Let $M$ be the $k$-marking of $Q:V_d=k$, $V_c=k_{2}$,
$V_b=k^{2}\sqcup k$ (see Example~1).

A representation of $Q_M$ may be
considered as a matrix over $k$ divided in the
following way:


\begin{center}
\begin{tabular}{c|c|c|c|c|c|}
\multicolumn{1}{c}{} & \multicolumn{2}{c|}{$c$} &
\multicolumn{2}{|c}{$b$} & \multicolumn{1}{c}{}\\
 \cline{2-6}
 $d$ & $c_1$ & $c_2$ & $x$ & $x^*$ & $y$\\
 \cline{2-6}
\multicolumn{1}{c}{} & \multicolumn{2}{c|}{$\rightarrow$} &
\multicolumn{2}{|c}{$\rightarrow $} & \multicolumn{1}{c}{}
\end{tabular}
\end{center}

Here $c_1>c_2$ correspond to the two objects of the partially
ordered set underlying $k_2$, so $\Hom(c_1,c_2)\neq 0$, and $x$
and $x^*$, with $xr=x^*$, correspond to the basis $\{x,x*\}$ of
the single object $X$ of $k^2$ with $\Hom (X,X) = k[r]/r^2$. The
strips $x$ and $x^*$ have the same numbers of columns, and when we
perform a column transformation inside of $x$, we should perform
the same column transformation in $x^*$. Moreover, we can perform
any row transformation, any column transformation inside $c_1$,
$c_2$, $y$, and add columns of $c_1$ to columns of $c_2$, and
columns of $x$ to columns of $x^*$.

So we see that representations of $k$-marked
quivers are obtained from representations of
non-marked quivers by division of matrices,
corresponding to the arrows, into vertical and
(in the general case) horizontal strips and the
restriction of admissible transformations. In
the case of non-$k$-marking, these matrices may
have fixed zero blocks.

Let us apply Corollary 2 to the arrow  $d \stackrel {\beta }
{\longrightarrow} c$, $Q'= d \stackrel {\gamma} {\longrightarrow}
b$, $\overline K_b=K_b=\oplus V_b$,  $ \overline K_d=\Rep (\Delta
_\beta )_M$. It is easy to see that $\ind (\Delta _\beta )_M$
consists of 5 representations $c^{0}_1$, $c^{0}_2$, $c^1_1$,
$c^1_2$, $d^0$, where $c^0_i\in k^{0\times 1}$, $c^1_i\in
k^{1\times 1}$, $d^0 \in k^{1\times 0}$ $(i=1,2)$. There arises a
functor $\Phi :\Rep (\Delta _{\beta})_M\longrightarrow \mod k$
such that $\Phi (c^0_i)=\left\{0\right\}$, and $\Im \Phi $ is
$\oplus{\overline{V}_d}$, $\Ob \overline{V}_d = \left\{c^1_1,
c^1_2, d_0 \right\}$, $\overline{V}_d \simeq k_3$.

Using Lemma 3 we conclude that there exists a bijection between
$\ind(Q_M)$ and $\ind (Q'_{\widetilde{M}}) \sqcup \left\{c^0_1,
c^0_2\right\};$ here $\widetilde{M}$ is a $k$-marking,
$\widetilde{V}(d) = k_3$, $\widetilde{V_b} = V_b$.
\end{example}
\citem A {\it bimodule}
$_{\mathfrak{A}}{\mathfrak{M}}_{\mathfrak{B}}$ over categories
$\mathfrak{A}$ and $\mathfrak{B}$ consists of
$\mathfrak{M}(A,B)\in \mod k$, $A\in \Ob\mathfrak{A}$, $B\in
\Ob\mathfrak{B}$ with natural multiplication and axioms. By $\El
\mathfrak{M}$ we denote a category whose objects are the elements
of the spaces $\mathfrak{M}(A,B)$, given $x\in \mathfrak{M}(A,B)$,
$y\in \mathfrak{M}(C,D)$, $\El \mathfrak{M} (x,y) = \{(\alpha,
\beta)/\alpha \in \mathfrak{A} (A,C)$, $\beta \in
\mathfrak{B}(B,D)$, $\alpha y = x \beta\}$. If $\mathfrak{A} =
\mathfrak{B}$ then $\widetilde{\El} \mathfrak{M}$ is the following
subcategory of $\El \mathfrak{M}: x\in \mathfrak{M}(A,B)$ belongs
to $\widetilde{\El} \mathfrak{M}$ if $A = B; (\alpha, \beta) \in
\Mor (\widetilde{\El} \mathfrak{M})$ if $\alpha = \beta$.

 Let $M$ be a $k$-marking of $Q$. Denote by
$D=D(Q_M)$ the Cartesian product of categories $K_x$, $x \in Q_v$.
Let $Q_v= \{q_1,\dots,q_n\}$; $X,Y \in \Ob D,
X=\{X_1,\dots,X_n\}$, $Y=\{Y_1,\dots,Y_n\}$, $X_i, Y_i \in \Ob
K_{q_i}$. If $q_i \stackrel {\alpha} {\longrightarrow}q_j$, put
${R(X,Y)}_\alpha=\mod k(X_i,Y_j)$. Let $R(X,Y)$ be the direct sum
of vector spaces $R(X,Y)_\alpha$, $\alpha \in Q_a$.

$R_\alpha (X,Y)$ may be naturally considered as a bimodule over
$K_x$, $K_y$ (if $_{x}\stackrel {\alpha} {\longrightarrow} _{y})$,
so bimodule $_{D}R_D$ arises.

It easy to see that the category $\widetilde{\El} {_{D}R_D}$ in
fact coincides with the category $\Rep Q_M$. Indeed an element
$U\in R(X,X)$ ``is'' the collection of $U(x) \in \Ob K_x$ $(x\in
Q_v)$ and $U(\alpha) \in \mod k$ $(U({x_\alpha}), U({y_\alpha}))$,
where $\alpha\in Q_a$, ${x_\alpha}
\stackrel{\alpha}{\longrightarrow} {y_\alpha}$, and morphisms in
$\Rep Q_M$ and $\widetilde{\El} {_{D}R_D}$ are defined in fact in
the same way.

So we may consider points (elements) of $R(X,X)$, $X\in \Ob D$ as
representations of $Q_M$ and talk about their indecomposibility
and equivalence ($\simeq$).

 A $k$-marked $Q_M$ is
 {\it wild} if for some $X$, $R(X,X)$
contains an affine plane $W_{X,X}$ consisting of indecomposable
and pairwise non-isomorphic representations.  A $k$-marked $Q_M$
is {\it tame} if $|\ind Q_M|=\infty$ and for any $X$ there exists
a finite set $T_{X,X}$ whose elements are affine subspaces of
dimension 0 and 1 in $R(X,X)$ such that each indecomposable
representation $U\in R(X,X)$ is equivalent to $w\in A\in
T_{X,X}$~[9].

 By $\aff k$ we denote the following $k$-category: $\Ob(\aff k) =
\Ob(\mod k)$, $\aff k(X,Y) = \{(A,y) |A \in \mod k(X,Y), y\in
Y\}$, $(A,y)(B,z) = (AB, yB+z)$. Given $x\in X$, $(A,y)\in \aff
k(X,Y)$ put $x(A,y) = xA+y$.

A functor $\Phi:\widetilde{\El}{_{\mathfrak{A}}
{\mathfrak{M}}_{\mathfrak{A}}}\longrightarrow
\widetilde{\El}{_{\mathfrak{B}}{\mathfrak{N}}_{\mathfrak{B}}}$ is
{\it affine} if there exists map $\overline\Phi: \Ob \mathfrak{A}
\longrightarrow \mathfrak{B}$ such that if $x\in \mathfrak{M}
(\mathfrak{A}, \mathfrak{A})$ then $\Phi(x)\in \mathfrak{N}
(\overline\Phi(\mathfrak{A}), \overline\Phi(\mathfrak{A}))$ and
the induced map $\Phi_A: {\mathfrak{M}}(A,A)\longrightarrow
{\mathfrak{N}}(\overline\Phi (A), \overline\Phi (A))$ is affine
($\in \Mor {\aff} k$) for each $A\in \Ob {\mathfrak{A}}$.

From our definitions directly follows

\begin{lemma} Let $_{D}R_D$ and  $_{D'}R'_{D'}$ be bimodules
attached respectively to $Q_M$ and $Q'_{M'},
\Phi:\widetilde{\El}{_{D'}{R'}_{D'}}\longrightarrow
\widetilde{\El}{_{D}R_{D}} $ is an affine functor. The wildness of
$Q'_{M'}$ implies the wildness of $Q_M$, if conditions 1) and 2)
hold:

1) $\Phi_X$ is embedding for each $X \in \Ob D'$;

2) $\Phi_X(x)$  is indecomposable if $x$ is, $x\simeq y$ if
$\Phi(x)\simeq \Phi(y)$ (for any $x,y\in R'(X,X)$, $X\in \Ob D'$).

The tameness of $Q'_{M'}$ implies the tameness of $Q_M$, if
conditions 3) and 4) hold:

3) for any $X\in \Ob D$, there exists only a finite number of
pairwise nonisomorphic objects $X'$ of $D'$, such that
$\overline{\Phi(X')}=X$;

4) if $X\in \Ob D$ then for any indecomposable $x\in R(X,X)$ exept
for a finite number there exists $x'\in R'(X',X')$ such that
$\Phi_{X'}(x')\simeq x$, $\overline{\Phi(X')} = X$.
\end{lemma}
Condition 2) clearly holds if, for each $X\in\Ob D'$, $x_1,x_2\in
R'(X,X)$ and $\varphi\in
\widetilde{\El}R(\Phi_X(x_1),\Phi_X(x_2))$, there exists
$\varphi{'}\in \widetilde{\El}R'(x_1,x_2)$, such that
$\Phi(\varphi{'}) = \varphi$.
\begin{remark} An affine functor
$\Phi:\widetilde{\El}{_{D'}{R'}_{D'}}\longrightarrow
\widetilde{\El}{_{D}R_{D}}$ naturally arises if $Q' = \Omega$, $M'
= \widetilde{\overline{M}}$ (see Proposition~2 and Lemma~3).
\end{remark}
It is easy to see that $Q_M$ can not be both wild and tame. In~[9]
it is proved that any vectroid is either wild or tame or has
finite type. In [10, 11] the same is proved for a wider class of
matrix problems, but in a formally different sense.

The marked quiver $Q_M$ of Example 5 is wild. Let $D\in \Ob K_d$,
$C\in \Ob K_c$ and
$B \in \Ob K_b$ all be 4-dimensional $k$-spaces. Then $R(X,X) =
\Hom_k(D,C) \oplus \Hom_k(D,B)$ contains the affine plane
\begin{displaymath}
W_{X,X}=
\begin{tabular}{|c|c|c|c|c|}
\multicolumn{2}{|c|}{c}
     &  \multicolumn{3}{|c|}{$b$}\\
     \multicolumn{5}{c}{ }\\[-4mm]
 0 & 0 & 1 0 & 0 0 & $\mu \lambda$ \\
 0 & 0 & 0 1 & 0 0 & 1 0 \\
 1 & 0 & 0 0 & 1 0 & 1 0 \\
 0 & 1 & 0 0 & 0 1 & 0 1\\
\multicolumn{2}{c}{$\rightarrow$}&
\multicolumn{2}{c}{$\rightarrow$}&
\multicolumn{1}{c}{}
\end{tabular}\ \ \
 \ \ (\lambda,\mu \in k),
\end{displaymath}
which consists of pairwise non-isomorphic indecomposable
representations of~$Q_M$.
\citem Since a vectroid is a subcategory
of $\mod k$, its objects are vector spaces, and its morphisms are
linear operators.

In [5], for any vectroid $V$, $\dim V= \sup \limits_{A\in \Ob
V}\dim A$ is defined. It is well known [3,~5] that if $V$ has
finite type then $\dim V\leq 3$. The rank of $V$, denoted $\rank
V$, is by definition a supremum of the ranks of noninvertible
additively indecomposable morphisms of $V$. Here a morphism
$\varphi $ is called {\it additively indecomposable} if $ \varphi
\neq \varphi _1 + \varphi _2$ where $\rank \varphi_1<\rank
\varphi$, $\rank \varphi_2<\rank \varphi$. It is known [5] that if
$\left|{\ind V}\right|<\infty$, then $\rank V \leq 2$. It is easy
to prove that if a vectroid $V$ is tame, then $\dim V\leq 4$,
$\rank V \leq 3$.

We define the {\it dimension}, $\dim Q_M$, of a $k$-marked quiver
$Q_M$ to be $\max\limits_{x\in Q_v}\dim V_x$.

It is clear that, if a $k$-marked quiver $Q_M$ has finite type,
then $\dim V_x \leq 3$ and $\rank V_x \leq 2$ for all $x\in Q_v$
and that, if it has tame type, then $\dim V_x \leq 4$ and $\rank
V_x \leq 3$ for all $x\in Q_v$.

It is easy to see that if $\dim V = 1$, then $V$ is a
linearization of some poset.

To an arbitrary vectroid $V$, we also associate a poset $S(V)$ [5]
in the following way. At first, consider the set $\overline{S}(V)$
consisting of all nonzero elements of all objects of $V$. Then we
define on $\overline{S}(V)$ the relation $<_\sim : \overline x
{<_\sim }\overline y$ if there exists $\varphi \in V(X,Y)$ such
that $\overline{x} \varphi = \overline{y} \quad (\overline{x} \in
X$, $\overline{y} \in Y$; $X, Y \in \Ob V)$. This relation is not
an ordering, because it can be that $x <_\sim y$ and $y <_\sim x$,
but after factorization by the induced equivalence we obtain a
poset $S(V)$.

We will also consider on $S(V)$ the equivalence relation $\simeq :
x \simeq  y$ if $\overline x \in X$, $\overline y \in Y$ and $X =
Y$, where $\overline x$, $\overline y$ are the elements of
$\overline {S}(V)$ corresponding to $x$, $y$; $X, Y \in \Ob V$. We
say that $x\in S(V)$ is {\it big}, if there exists $y\simeq x$,
$y\neq x$.

We will say that the vectroid $V$ is {\it halflinear}
if it is not linear and

1)   $\dim V \leq 2$

2)   $\rank V = 1$

3)   if $a$ is big then $a$ is comparable with
each point of $S(V)$

4)   if $a$ is small then $a$ can be
incomparable with only one point of $S(V)$.

It is easy to see that a halflinear (or linear) vectroid $V$ is
determined uniquely by $(S(V), \leq , \simeq )$.

A vertex $x \in Q_v$ is said to be {\it linearly marked} if
$V_x$ is linear and to be {\it halflinearly marked} if $V_x$
is halflinear. The marked quiver $Q_M$ is said to be
{\it linear} if each $x\in Q_v$ is linearly marked, and to
be half linear if each $x$ is either halflinearly or linearly
marked.

By $G(Q)$ we  denote the non-oriented graph that corresponds to a
quiver~$Q$; we will view a $k$-marking of $Q$ also as a marking of
$G(Q)$ -- associating to each $x\in G(Q)$ a vectroid $V_x$. If all
$V_x$ are linear or halflinear, we will write over each vertex $x$
of $G(Q)$ the number $n$ if $V_x = k_n$, the symbol $\infty$ if
$V_x$ is halflinear, and nothing  if $V_x = k$.

A $k$-marked quiver $Q_M$ is {\it Gelfand}
 if and only if either $Q = \Delta $, and both vectroids
$V_x, V_y$ are halflinear, or $G(Q)_M = \bullet$
--- $\overset{\infty} \bullet$ --- $\bullet$.

Representations of the Gelfand $k$-marked quivers were treated  in
[12], where in fact their tameness was proved. (If $Q = z
\longleftarrow x \longrightarrow y$, then $\Rep Q_M = \Rep
\Delta_{\overline{M}}$, where $\overline{V}_x = V_x$,
$\overline{V_y} = k \sqcup k$. If $Q =z\longrightarrow x
\longrightarrow y$, then Corollary 2 implies $\Rep Q_M = \Rep (z
\longleftarrow x \longrightarrow y)_{ \overline{M}}$ where
$\overline{V_{x}} = V_{x}^{\circ}$, $(\overline{V_{y}} =
\overline{V_{z}} = k)$.

For a halflinear vectroid $V$, we denote by $V^-$ the vectroid
obtained from~$V$ by excluding those one-dimensional objects which
correspond to points of $S(V)$ comparable to all other points of
$S(V)$. We say that two halflinear vectroids $V_1$ and $V_2$ are
{\it almost equivalent} if $V^-_1 \simeq V^-_2$.
\begin{lemma} $Q_M$ is wild if it is a $k$-marked quiver for
which $G(Q)_M$ is one of the the following types:

1) $\bullet$\raisebox{1mm}{\rule{5mm}{0.4pt}}$
\underset{x} \bullet$\raisebox{1mm}{\rule{5mm}{0.4pt}}$\underset{y}
\bullet$, where $V_x$ is not linear, $V_y\neq k$.

2) $\bullet$\raisebox{1mm}{\rule{5mm}{0.4pt}}$\underset{x} \bullet$
\raisebox{1mm}{\rule{5mm}{0.4pt}}$\bullet$, where
$V_x$  is not linear and not halflinear.

3) $\underset{x}
\bullet$\raisebox{1mm}{\rule{5mm}{0.4pt}}$\underset{y} \bullet$,
where both $V_x$, $V_y$  are not linear, and $V_x$ is not
halflinear.

4) a cycle $\tilde{A}_n$ containing a vertex $a$ with $V_a\neq k$.
\end{lemma}
The proof is straightforward.
\begin{remark}  Let $Q_M$ be a $k$-marked quiver with an arrow $x
\stackrel{\alpha}{\longrightarrow} y$, $V_x = V_y$ and let $L_N$
be the $k$-marked quiver obtained from $Q_M$ by excluding the
arrow $\alpha$ and uniting the vertices $x$ and $y$ in one vertex.
Then $\Rep Q_M$ contains a full subcategory $R' \simeq \Rep L_N$.
Indeed, consider the full subcategory $R' \subset \Rep Q_M$,
defined by $U\subseteq \Ob R'$ if $U(x) = U(y)$, $U(\alpha) =
1_{U_x}$. So if $L_N$ is wild, then $Q_M$ also is.
\end{remark}
\citem Let $\beta$ be the arrow of  the $k$-marked $Q_M$
considered in Section~2, that is, $Q$ contains a vertex $z$ with
$\delta(z) =1$ and either $z\stackrel{\beta}{\longrightarrow} w$
or $w \stackrel{\beta}{\longrightarrow} z$, with $w\neq z$ and
$|Q_a|>1$. We will say that $\beta$ is {\it reducible} in one of
the following cases:

1) $V_w = k$, $V_z = k_m$

2) $V_w = k$, $V_z$ is halflinear

3) $V_z = k$, $V_w = k_2$.
\begin{lemma} If $\beta$ is reducible, then there exists a natural
bijection between $\ind Q_M \setminus \Ob V_z$ and $\ind Q'_{M'}$.
If $Q'_{M'}$ is tame, then so is $Q_M$; if it is wild then $Q_M$ is wild.
\end{lemma}
Here $Q'$ is as defined in Section~2. Thus $V'(x) = V(x)$ if $x\neq
w$. Furthermore $V'(w) = k_{m+1}$ in case 1); $V'_w$ is halflinear
and $V'_w$ is almost equivalent to $V_w$ in case 2); $V_w$ is
halflinear and $S(V_w) = kP$, where P is the poset $\{p_1, p_2,
p_3, p_4 | p_1 < p_2 < p_4, p_1 < p_3 < p_4 \}$ in case~3).

The proof follows from Corollary~2, Lemma~3, Lemma~5, Remark~5 and
the calculation of $\Rep (\Delta_\beta)_M$ in cases 1.2.3.
\begin{example} Let $Q$ be a quiver containing subquiver
$Q^0= \Delta=w \longrightarrow z$, such that $\delta(z)=1$, and
$M$ be such marking of $Q$ that we have the situation of the
Lemma, case~2). Namely $K_w=\mod k$, $V_z$ is halflinear, $\Ob
V_z=\{A,B,C,D\}$, $\dim A=2$, $\dim B=\dim C=\dim D=1$,
$S(V_z)=\{a,a^*,b,c,d \}$, $a<b<a^*<c$, $a^*<d$. Let
$Y^0_w:\Rep\Delta \longrightarrow K_w=\mod k$ be the functor
determined in Section~2. We may consider $\Ob V_z\subset \ind
\Delta_M$. Let $T\in \ind \Delta_M$, then $Y^0_w(T)=0$ if and only
if $T\in \Ob V_z$. Here $\ind \Delta \backslash \Ob
V_z=\{a_1,a^*_1, \overline A, (c,d), b_1, c_1, d_1, w_0 \}$.

$a_1 = (\underset{a} 1 | \underset{a^*} 0)$, $a^{*}_1 =
(\underset{a} 0 | \underset{a^*} 1)$,

\begin{displaymath}
\mbox{\raisebox{2mm}{$\overline{A}={}$}}
\begin{tabular}{|c|c|}
$1 \vert 0$\\ $0 \vert 1$ \\ \multicolumn{2}{c}{$a$  $a^*$}
\end{tabular},
\end{displaymath}
$(c, d) = (\underset{c} 1 | \underset{d} 1)$, $t_1 = (\underset{t}
1)$ for $t \in \{b, c, d\}$.

After application of Corollary 2 and Lemma 3 we get $Q'_{M'}$
where $\Ob V'_w=Y^0_w (\ind \Delta\backslash \Ob V_z)=\{w_0, b_1,
c_1, d_1, a_1, a^*_1, \overline A, (c,d)\}$ where $\dim \overline
A = 2$, (the dimensions of the others are equal to 1). Furthermore
the poset $S(V'_w) = \{w_0, b_1, c_1, d_1, a_1, a^*_1, \overline
a, \overline a^*, (c, d)\}$, $(\overline a < \overline a^*)$, $a_1
< \overline a <b_1 < a^*_1 < \overline a^* < (c,d) < d_1 < w_0$,
$c, d<c_1,w_0$. The vectroid $V'_w$ is clearly almost equivalent
to $V_w$. It is easy to see that the conditions of Lemma~5 (in
this and the general cases) hold.
\end{example}
\begin{proposition}
Let $Q_M$ be $k$-marked, $Q = \Delta_\alpha$, ${x
\stackrel{\alpha}{\longrightarrow}y}$ and $V_x = k_m$, $m > 1$.
Then there exists a natural bijection between $\ind Q_M$ and $\ind
V\setminus \Ob k_{m- 1}$, where $V = V_y\sqcup k_{m-1}$. If $Q_M$
is tame, $V$ is tame; if $Q_M$ is wild, $V$ is wild. If $V$ is
tame then $Q_M$ also is.
\end{proposition}
\begin{proof}
We consider $\overline Q = y \longleftarrow x \longrightarrow z
\quad \overline V_x = k$, $\overline V_y = V_y$, $\overline V_z =
k_{m-1}$. Then the first and the second statements follow from the
application of the Lemma to the arrow $x\longrightarrow z$.

In order to prove that tameness of $V$ implies tameness of $Q_M$
we use a matrix construction, although it probably should follow
from some general consideration. In the picture below we will for
simplicity assume  $m=3$. The representation $U$ of $V_y\sqcup
k_2$ has the form
\begin{center}
\begin{tabular}{|ccc|c|c|}
\hline & $U'$ & & ~ & ~
 \\ \hline
\multicolumn{5}{c}{$\hspace{17mm}\longrightarrow $} \\
\end{tabular}
\end{center}
where $U'$ is a representation of $V_y$ and the empty columns are
for a representation of $k_2$ which we reduce to a standard form
in the next diagram.

We say that a representation $U$ has
{\it preliminary} form if

\begin{displaymath}
\mbox{\raisebox{2mm}{$U={}$}}
\begin{tabular}{|c|c|c|}
\cline{1-3}
$ U_1$ & $U_{11}$& $U_{12}$ \\ \cline{1-3}
$U_2$ & 0 & $U_{22}$\\ \cline{1-3}
${U_3}$& 0 & 0\\
\cline{1-3}
\multicolumn{1}{c} {} &
 \multicolumn{2}{c} {$\longrightarrow$}
\end{tabular}
\end{displaymath}
where $U_{11}$ and $U_{22}$ are matrices with
linearly independent rows. Let $r_t(U)$ be the number of rows of
 $U_t$ $(t=1,2,3)$.
Of course any representation may be reduced to preliminary form.

Let $P(V)$ be a full subcategory of $\Rep V$ consisting of
representations in preliminary form.

Consider a map $F: \Ob P(V) \longrightarrow \Ob \Rep Q_M$:
\begin{displaymath}
F(U)=
\begin{tabular}{|c|}
$ U_1$  \\
\hline
$U_2$ \\
\hline
${U_3}$
\end{tabular}
\begin{array}{c}\uparrow \\[1mm] \uparrow \end{array}
\end{displaymath}

It is clear that for each $Z'\in \Ob D(Q_M)$ (see Section~5) there
exists $Z\in \Ob D(V)$ such that $F(R(Z,Z)\bigcap P(V))=
R(Z',Z')$. (If $Q=\Delta$, $Z'=(X,Y)\in \Ob D(Q_M)$, $X\in K_x$,
$Y\in K_y$ then $R(Z',Z')\simeq \Hom_k(X,Y)$).

It is easy to see (using Corrolary~2 or directly) that $F$ can be
considerred as a functor.

Let $A$ be an affine subspace of $R(Z,Z)$, $Z\in \Ob D(V)$, $\dim
A=1$. We will write $A\widetilde\subset P(V)$ if
$A\supset\overline A$, $|A\backslash \overline A|< \infty$,
$\overline A \subset P(V)$ and $r_t(a_i)=r_t(a_j)$ for $a_i,
a_j\in \overline A$ $(t=1,2,3)$. It is clear that if
$A\widetilde\subset P(V)$ then there exists an affine subspace
$\overline {F(A)}$ of $R(Z',Z')$ such that $\dim \overline
{F(A)}\leq 1$, $|\overline {F(A)}\backslash F(\overline
A)|<\infty$.

It is easy to show that row-transformations may be used to reduce
any $A$ to (an affine subspace) $A'\widetilde\subset P(V)$ (i.e.
an ivertable matrix $M$ exists such that $MA=A'$). Now we may
consider that if $V$ is tame, then for each $Z\in \Ob D(V)$ there
exist $A_1,\dots,A_n \widetilde\subset P(V)$ $(\dim A_i=1)$ which
generate almost all indecomposable points of $R(Z,Z)$.
Consequently $\overline{F(A_1)},\dots,\overline{F(A_n)}$ generate
almost indecomposable points of all those $R(Z',Z')$ for which
$F(R(Z,Z)\cap P(V))=R(Z',Z')$. So $Q_M$ is tame with $V$.
\end{proof}
\begin{remark} A similar method may be used to prove that, in the
situation of the Lemma, case~1, the tameness of $Q_M$ implies the
tameness of $Q'_{M'}$.
\end{remark}
For a halflinear $k$-marked $Q_M$ we construct a (non-oriented)
graph $G(Q_M)$ in the following way. For each vertex $x$ such that
$V_x = k_n$, we add to the graph $G(Q)$ $(n-1)$ vertices $a^x_2
\cdot\cdot\cdot a^x_n$, $x$ --- $a^x_2$ --- $\cdots$
--- $a^x_n$,  for each vertex $y$ such that $V_y$ is halflinear we add to
$G(Q)$ two vertices $ b^y_1$, $b^y_2$, $b^y_1$ --- $y$
--- $b^y_2$.

It is easy to see that under the conditions of the Lemma in cases
1 and~3, $G(Q_M)\simeq G(Q'_{M'})$. In case 2, if $G(Q_M)$ is a
Dynkin or extended Dynkin diagram, then $G(Q_M)$ is $D_n$ or
$\tilde{D}_n$, with $n>4$, respectively, and $G(Q'_{M'})$ is
$D_{n-1}$ or $\tilde{D}_{n-1}$, respectively. Conversely, if
$G(Q'_{M'})$ is a Dynkin or extended Dynkin diagram, then
$G(Q'_{M'})$ is $D_n$ or $\tilde{D}_n$, with $n>3$, and $G(Q_M)$
is $D_{n+1}$ or $\tilde{D}_{n+1}$, respectively.
\begin{theorem}
If $Q_M$ is halflinear, then it is tame, or has finite type, if
and only if $G(Q_M)$ is an extended Dynkin diagram, or a Dynkin
diagram, respectively.

Suppose there is a vertex $a\in Q_M$ with $V_a$ neither a
halflinear nor a linear vectroid. Then $Q_M$ is tame, or has
finite type, if and only if $G(Q) = A_n$ with $a=a_1$, $V_{a_i}
=k$ for $1<i<n$ and $V_n= k_m$ for some $m\geq 1$, and the
vectroid $V\sqcup k_{m+n-3}$ is tame, or has finite type,
respectively; here $V= V_a$ or $V_a^\circ$ according as $a$ is the
end or beginning of an arrow of $Q_M$.
\end{theorem}
\begin{proof}
We will assume that for unmarked quivers the statement is known.
Although usually tameness is defined in a different way, it is
easy to see that the proofs are valid for tameness in our sense.

So we may consider that $G(Q)$ is a Dynkin diagram or an extended
Dynkin diagram.

We will also assume that $G(Q)$ (and so $G(Q_M)$) is acyclic
because an unmarked cycle is tame and a marked cycle $Q_M$ is wild
by Lemma~6 (case~4) (and $G(Q_M)$ is neither a Dynkin nor an
extended Dynkin diagram).

a) Let $Q_M$ be tame, or of finite type $a\in Q_v$, $V_a$ be
neither linear nor halflinear. Then $\delta(a)=1$ by Lemma~6~(2).

$a_1$) $Q_v \ni x \neq a$, $V_x$ is not linear. Using Remark~6 we
have a contradiction to Lemma 6 (3).

$a_2$) $Q_v \ni x \neq a$, $V_x \neq k$, $\delta (x) \neq 1$. We
show that $Q_M$ is wild (that contradicts which our assumption).

Without loosing of the generality we can (using Remark 6) assume
that \\ $G(Q)=$  $\overset{a} \bullet$
--- $\overset{x} \bullet$ ---  $\overset{y} \bullet$ , $V_x=k_2$, $V_y=k$.

Using the Lemma (case 3) we get $Q'_{M'}= $  $\overset{a} \bullet$
--- $\overset{x} \bullet$,
where $V'_x$ is halflinear ($V'_a=V_a$), $Q'_{M'}$ is wild by
Lemma 6 (3) and so also is $Q_M$ (the Lemma).

$a_3$) $Q_v\ni x, \delta(x)>2$. This case is reduced to $a_2$) by
several applications of the Lemma (case 1) and Remark~6.

$a_4$) $Q=A_n$, $a=a_1$, $V_{a_2}=\cdots=V_{a_{n-1}}=k$,
$V_{a_n}=k_m$. By $n-2$ applications of the Lemma (case 1) and (at
the end) of the proposition reduce $\Rep Q_M$ to $\Rep (V_a\sqcup
k_{n+m-3})$ or $\Rep (V^{\circ}_a \sqcup k_{n+m-3})$.

Now we assume that $Q_M$ is halflinear.

b) Let $Q_v \ni b, V_b$ be halflinear.

$b_1$) $\delta (b)\ne 1$. If $G(Q)= $ $ \bullet$
--- $\overset{\infty} \bullet$ --- $\bullet$ ($Q_M$ is
Gelfand, $G(Q_M)=\Tilde{D}_4$), then it is tame. (If $Q$ is $
\bullet \longrightarrow
 \bullet \longleftarrow \bullet $ or $
\bullet \longleftarrow \bullet \longrightarrow \bullet $ then
$\Rep Q_M \simeq \Rep \Delta_{\overline {M}}$ where
$\overline{V}_x = V_b, \overline{V}_y = k\sqcup k$ or
$\overline{V}_x = k\sqcup k, \overline{V}_y = V_b$, so
$\Delta_{\overline{M}} $ is tame by the results of~[12]. If $Q =
\bullet \longrightarrow
 \bullet \longrightarrow \bullet $, then $\Rep Q_M \simeq \Rep Q^{*}_M$
 where $Q^{*}=  \overset{x}
\bullet \longrightarrow \overset{b} \bullet \longleftarrow
\overset{y} \bullet$, $V_x = V_y = k$, $V^{*}_b = V^{\circ}_b)$.

Suppose $Q_M$ is not Gelfand. Then it is easy to see that $G(Q_M)$
is not an extended Dynkin (and not Dynkin) diagram. If $Q_v\ni x$,
$V_x \neq k$ then $Q_M$ is wild by Lemma~6~(1) and Remark~6. If
all vertices except for $b$ are unmarked but $|Q_v|>3$, we reduce
this case to the case above using the Lemma (case~1).

$b_2$) $\delta (b) = 1$. In this case we will prove the theorem by
induction on $|Q_v|$. If $|Q_v| = 2$, $Q_M$ is Gelfand and tame
(see~[12]) or of finite type, then $G(Q_M)$ is $\Tilde{D}_5$ or
$D_n$ respectively. So let $|Q_v|>2$ and suppose that
$$ G(Q_M)
\supset\overset{b}\bullet \mbox{---} \overset{x} \bullet  \mbox{---}
\overset{y} \bullet. $$
 If $V_x = k$ we apply the Lemma (case~2). Now $G(Q'_{M'})$ is an
extended
 Dynkin (Dynkin) diagram if and only if
$G(Q_M)$ is. Since $|Q'_v|<|Q_v|$ our statement follows from the
inductive assumption.

Let $V_x = k_m$, $m>1 (x$ can not be halflinear by $b_1$). If $m
\geq 3$ then $G(Q_M)$ is neither an extended Dynkin nor a Dynkin
diagram and we show that $Q_M$ is not tame. We may assume $Q_M = $
$\overset{\infty} \bullet$
--- $\overset{3} \bullet$ ---  $\bullet$.

Consider $L_N =$  $\overset{\infty}{b} $
--- $
\overset{\overset{\displaystyle \overset{2} z }|} w $ --- $x$.
Then $Q_M = L'_{N'}$ (the Lemma, case 1).

On the other hand we apply the Lemma, case 2 to the vertex $b$ of
$L_N$ and get  $\overset{2} \bullet$
--- $\overset{\infty} \bullet$ ---  $ \bullet$. The last
$k$-marked quiver is wild by Lemma~6~(1), so also is $L_N$ (by the
Lemma, case 2). Now $Q_M$ is not tame because otherwise $L_N$
should be tame by the Lemma (case~1). Let $m=2$. The cases $V_y
\neq k$ and $|Q_v|>3$ are considered on the analogy of the cases
$m\geq 3$ and $b_1$).

If $G(Q_M) =$  $\overset{\infty} \bullet$
--- $\overset{2} \bullet$ ---  $ \bullet$ we get a tame
$Q_M$ by the Lemma, case 3 and~[12] ($G(Q_M) = \widetilde{D}_5$).

c) $Q_M$ is  linear. We consider
an unmarked quiver $\overline{Q}$ such that
$G(\overline{Q}) = G(Q_M)$ and by the several
applications of the Lemma  case~1 to $\overline {Q}$ we reduce it
to $Q_M$. The remark implies than if $\overline{Q}$ is tame or
of finite type, then $Q_M$ is of the same type.  The converse statement
follows from the Lemma.
\end{proof}
The theorem gives a criterion for wildeness and finiteness of type
for halflinear (in particular linear) $Q_M$, and for a $k$-marked
quiver of dimension~1 (together with [13, 14]), as well as a
criterion for the finiteness of type if $\dim Q_M = 2$ (together
with [15, 8], see also [3, 16--20]).

If $\dim Q_M \geq 3$ and $|\ind Q_M| < \infty$, then $\dim Q_M =
3$ (see Section~6), and it is easy to see that in this case $Q =
\Delta$, and if $Q_v = \{x,y\}$ and $\dim V_y = 3$, then $V_x =
k$, so we have representations of a vectroid $V$, $\dim V = 3$.
For this case a criterion for the finiteness of type is formulated
in~[21], but only the necessity of it has been proven.

  From our considerations in
fact it follows that the tameness in our sense coincides with the
tameness in sense of [10, 11]. In fact Lemma~5, Lemma~7 and
Proposition~7 hold for tameness and wildness in the sense of
[10, 11] (the last statement in the Proposition~7 follows from the
second by Drozd's theorem).
\end{pgraph}
\smallskip

{\large \bf Appendix}
\smallskip

Here I attempt to present my personal point of
view on the trends in investigations connected
with the subject of this article.

In the first half of the XXth century it was the opinion of a
majority of mathematicians that after the classical results of
Jordan, Weierstrass and Kronecker, linear algebra as a branch of
pure mathematics was completely finished. As an exception, one
should mention the remarkable work of Szekerez~[22].

Concerning the classical results mentioned above, it seemed that
in spite of their importance for applications they were far from
abstract algebra, and in particular from representation theory. I
remember that 40 years ago all Soviet algebraists were surprised
and even in some sense disappointed when V.A. Bashev (a student of
I.R. Shafarevich) had solved~[23] the problem of classification of
representations of the Klein 4-group over a field of
characteristic 2 (which was considered then as very difficult) by
a trivial reduction of this problem to the problem of Kronecker.

Between 1960 and 1970 several different problems were solved by
their reduction to some classification problems of linear algebra
[24--26].

The concept of matrix problems (or combinatorial
problems of linear algebra) as a special branch
of mathematics arose about 1970
from at least three sources.

In [27, 28] (see also [29]) it was conjectured in particular that
the category of Harish-Chandra modules for any semisimple group is
equivalent to a certain category of diagrams in the category of
finite dimensional vector spaces, and a boom in linear algebra was
predicted on the basis of the new categorical and homological
methods.

In [30, 31, 2] representations of posets were introduced, and it
was claimed that many other matrix problems can and should be
reduced to them.

At last, but of course not least, it was clear
that the subject of [1] was very wide and deep.
I hope that this article of mine underlines once more
the importance of representations of quivers.

In the coming years the theories of representations of quivers and
posets were developed successfully.  It was clear from the
beginning  that these two theories were very close.
Representations of a majority of quivers can be reduced to
representations of posets (see for example [32]). Note that ``on
the way'' from representations of quivers to representations of
posets, in fact, representations of marked quivers arose. Many
important generalizations of representations of posets also were
introduced. In particular, representations of vectroids played a
big role in many questions. P. Gabriel showed that representations
of arbitrary finite dimensional algebra can be reduced to
them~[9].

However the general theory of matrix problems was not developed so
well. It may be partially explained by the absence of a natural
basic definition. The widest class of matrix problems are the
representations of $\mathcal{DGC}$ or bocses, which were
introduced by M. Kleiner and me [33]. But this class is not only
wide enough, but may be in some sense too wide. It includes the
representations of posets, vectroids, and ``Gelfand problems'',
but these most important problems are not picked out naturally in
the very wide area of bocses. The general definition of bocses is
more or less clear (though is not so natural as quivers and
posets), but additional conditions which are necessary to get a
majority of really important classes of matrix problems are rather
complicated. So the theory of bocses became convenient to prove
general theorems [10, 11], but not to develop a systematic theory
of different classes of matrix problems.

Conversely, in the terms of this article the
most important classes of matrix problems are
picked out very naturally, and it seems that the
general definition can stimulate investigations
of many other natural and useful classes of
matrix problems. So it seems that the notion of
representations of marked quivers should be a
better basic notion for classification problems
of linear algebra.

I want to underline that in spite of the close
relation between matrix problems and
representations of algebras, it does not seems
correct to consider the first theory only as
part of the second. Many matrix problems arose
from applications to other branches of
mathematics, and we hope that the number of such
applications will increase in the future.

The first aim in the theory of matrix problems is to finish a
description of vectroids of finite type and their representations.
As follows from this article, this will imply that such a theory
also exists for representations of $k$-marked quivers. For the
point-marked and non point-marked (see Example~4) quivers such
theories should be more varied but also solid.

In the theory of representations of tame vectroids (and tame
$k$-mar\-ked quivers) the first important class is
{\it locally semisimple} vectroids (i.e. $\Hom(A,A)$ is a
semisimple algebra for any $A$)\footnote[1]{Such categories were
already considered by J.M. Gelfand and G.E. Shilow in 1963 [34].}. A
criterion of tameness for such vectroids was announced by
L.A. Nasarova in 1985 (see [35]).

At the end I want to note that although a
majority of matrix problems are wild, their
investigations may make sense. For any matrix
problem (representations of quiver, poset,
marked quiver, bocs) all representations of
fixed dimension form a vector-space. Any affine
subspace of this space generates some sets of
representations of this and bigger dimensions.

\noindent
{\bf Conjecture.}  {\it For any matrix problem there exists an
$n$ such that all indecomposable representations of fixed
dimension are generated by a finite number of affine subspaces of
dimensions at most $n$, whose elements are indecomposable and
pairwise inequivalent.}

For $n = 0$ it would imply the first Brauer-Thrall conjecture
[36], and for $n = 1$, -- Drozd's theorem  [10, 11, 9].

Of course this formulation may be too strong,
but it seems that some its modifications will be
typical for the investigations in the new
millenium.

\newpage
\centerline{\large \bf References}

\begin{enumerate}
\itemsep0mm
\item Gabriel P.: Unzerlegbare Darstellungen I,
{\it Manuscr. Math.} {\bf 6} (1972), 71--103.

\item Nazarova L.A., Roiter A.V.: Representations of partially ordered
sets, {\it Zap. Nauchn. Semin. Leningr. Otd. Mat. Inst. Steklova}
{\bf 28} (1972), 5--31, English transl.: {\it J. Sov. Math.} {\bf
3} (1975), 585--606.

\item Gabriel P., Roiter A.V.: Representations of Finite-Dimensional
Algebras: Springer-Verlag, Algebra VIII, 1992.

\item Bass H.: Algebraic K-theory. W.A. Benjamin, INC. New York,
1969.

\item Belousov K.I., Nazarova L.A, Roiter A.V., Sergeichuk V.V.:
Elementary and multielementary representations of vectroids, {\it
Ukr. Math. J.} {\bf 47} (1995), 1451--1477.

\item Nazarova L.A., Roiter A.V.: Kategorielle Matrizen-Problems und
Brauer-Thrall - Vermutung,  Mitt. Aus dem. Math. Sem. Giessen,
1975.

\item Ringel C.M.: Tame algebras and integral quadratic form,
{\it  Lecture Notes Math.} {\bf 1099} (1984).

\item Nazarova L.A., Roiter A.V.: Finitely represented dyadic sets,
{\it Ukr. Math. J.} {\bf 52} (2000), 1363--1396.

\item Gabriel P., Nazarova L.A., Roiter A.V., Sergeichuk V.V., Vossieck
V.: Tame and wild subspace problems, {\it Ukr. Math. J.} {\bf 45}
(1993), 313--352.

\item Drozd Ju.A.: Tame and wild matrix problems (Ottava, 1979),
{\it Lect. Notes Math.} {\bf 832} (1980), 242--258.

\item Crawley - Boevey W.W.: On tame algebras and bocses,
{\it Proc. Lond. Math. Soc. III Ser.} {\bf 56} (1988), 451--483.

\item Nazarova L.A., Roiter A.V.: On a problem of I.M. Gelfand,
{\it Funk. Anal. and Appl.}, {\bf 7} (1973), 54--69.

\item Nazarova L.A.: Partially ordered sets of infinite type,
{\it Izv. A.C. USSR}, {\bf 39} (1975), 963--991.

\item Kleiner M.M.: Partially ordered sets of finite type,
{\it Zap. Nauchn. Semin. Leningr. Otd. Math. Inst. Steclova}, {\bf
281} (1972), 32--41.

\item Roiter A.V., Belousov K.I., Nazarova L.A.:
Representations of finitely represented dyadic sets, {\it Algebras
and modules II, ICRA 1996 CMS Conf. Proc.} {\bf 24} (1996),
61--77.

\item Nazarova L.A., Roiter A.V.: Representations of biinvolutive
posetsI, Preprint 91.34, Kiev, Math. Inst., 1991.

\item Nazarova L.A., Roiter A.V.: Representations of biinvolutiv posets
II, Preprint, Kiev, Math Inst., 1994.

\item Guidon T.: Representations of dyadic sets 2, Diss. Uni.
Zurich, 1996, 1--47.

\item Hassler U: Representations of dyadic sets 1, Diss. Uni.
Zurich, 1996, 1--62.

\item Guidon T., Hassler U., Nazarova L.A., Roiter A.V.: Dyadic sets S:
a dichotomy for indecomposable S-matrices, Comptes - Rendus Acad.
Sc. Paris, {\bf 324}, Series I (1997), 1205--1210.

\item Belousov K.I., Nazarova L.A., Roiter A.V.: Finitely represented
tryadic sets, {\it Alg. and Analis} {\bf 9(4)} (1997), 3--27.

\item Szekeres G.: Determination of certain family of finite metabelian
groups, {\it Trans. Amer. Math. Soc.} {\bf 66} (1949), 11--43.

\item Bashev V.A.: Representations of group $Z_2\times Z_2$ in a field
of characteristic 2, {\it DAN SSSR} {\bf 141:5} (1961),
1015--1018.

\item Nazarova L.A.: Integral representations of Klein's four group,
{\it Dokl. Akad. Nauk SSSR} {\bf 140} (1961), 1011--1014;
Representations of the local ring of a curve with 4 branches, {\it
Izv. Akad. Nauk SSSR. Ser. Math.} {\bf 31} (1967), 1361--1378.

\item Drozd Yu.A., Roiter A.V.: Commutative rings with a finite number
of integral indecomposable representations, {\it Izv. Akad. Nauk
SSSR} {\bf 31} (1967), 783--798.

\item Gelfand I.M., Ponomarev V.A.: Indecomposable representations of
Lorentz group, {\it Usp. Math. Nauk} {\bf 32} (1968), 3--60.

\item Gelfand I.M.: The cohomology of infinite dimensional Lie
algebras, some questions of integral geometry, {\it  Actes I.C.
Math.} {\bf 1/10} (1970) /NICE/ France 95--111.

\item Gelfand I.M., Ponomarev V.A.: Problems of linear algebras and
classification of quadruples of subspaces in finite-dimensional
vector space, Colloq. math. Soc. Janos Bolyai 5 (Tihany 1970)
(1972) 163--237;  {\it DAN SSSR} {\bf  197:4} (1971), 762--765.

\item Bernstein I.N., Gelfand I.M., Ponomarev V.A.: Coxeter functors
and Gabriel's theorem, {\it Izv. Mat. Nauk} {\bf 28} (1973),
19--33; English transl. {\it  Russ. Math. Surv.} {\bf 28} (1973),
17--32.

\item Nazarova L.A., Roiter A.V.: Matrix Questions and the
Brauer-Thrall Conjectures on Algebras with an Infinite Number of
Indecomposable representations, {\it Proc. Amer. Math. Soc.}
(1971), 111--115.

\item Roiter A.V.: Matrix problems and the representations of
bisystems, {\it Zap. Nauchn. Sem. LOMI} {\bf 29} (1972), 130--139.

\item  Nazarova L.A.: Representations of quivers of
infinite type, {\it Izv. AN SSSR} {\bf 37:4} (1973), 752--791.

\item Kleiner M.M., Roiter A.V.: Representations of D.G.C.,
in Matrix problems,  Kiev, Inst. Math. AN USSR, 1977, 5--71

\item Gelfand I.M., Shilov G.E.: Categories of finite-dimensional
spaces, {\it Usp. Math. Nauk Ser.} {\bf 1} (1963), 27--49.

\item Roiter A.V.: Representations of posets and tame matrix problems,
in Proc. Durham Simp. 1985, Cambridge, 1986, 91--109.

\item Roiter A.V.: Unboundness of the dimension of the indecomposable
representations of algebras, that have an infinite number of
indecomposable representations, {\it Izv. AN SSSR} {\bf 32}
(1968), 1275--1282.
\end{enumerate}

\end{document}